\documentclass[12pt]{amsart}

\usepackage{amsfonts,amsmath,amsthm,amssymb}
\usepackage{latexsym}

\newtheorem{theorem}{Theorem}[section]

\newtheorem{proposition}[theorem]{Proposition}

\theoremstyle{definition}

\newcommand{\Z}{\ensuremath{\mathbb{Z}}}
\newcommand{\R}{\ensuremath{\mathbb{R}}}

\newcommand{\T}{\ensuremath{\mathcal{T}}}

\def \C {\Gamma}

\def \f {\phi}

\def \< {\langle}
\def \> {\rangle}

\newcommand{\beql}[1]{\begin{equation}\label{#1}}
\newcommand{\eeq}{\end{equation}}
\newcommand{\comment}[1]{}

\newcommand{\Abs}[1]{{\left|{#1}\right|}}

\newcommand{\Norm}[1]{{\left\|{#1}\right\|}}

\newcommand{\Set}[1]{{\left\{{#1}\right\}}}

\newcommand{\RR}{{\mathbb R}}

\newcommand{\ZZ}{{\mathbb Z}}

\newcommand{\TT}{{\mathbb T}}

\newcommand{\inner}[2]{{\langle #1, #2 \rangle}}

\newcommand{\ft}[1]{\widehat{#1}}

\renewcommand{\vec}[1]{{\mbox{\boldmath$#1$}}}

\newcounter{rem}
\begin{document}

\title{Tiles with no spectra}

\author{Mihail N. Kolountzakis \& M\'at\'e Matolcsi}

\date{June 2004}

\address{M.K.: Department of Mathematics, University of Crete, Knossos Ave.,
GR-714 09, Iraklio, Greece} \email{kolount@member.ams.org}

\address{M.M.: Alfr\'ed R\'enyi Institute of Mathematics,
Hungarian Academy of Sciences POB 127 H-1364 Budapest, Hungary.}\email{matomate@renyi.hu}

\thanks{Partially supported by European Commission IHP Network HARP
(Harmonic Analysis and Related Problems), Contract Number: HPRN-CT-2001-00273 - HARP}

\begin{abstract}
We exhibit a subset of a finite Abelian group, which tiles the group by translation, and
such that its tiling complements do not have a common spectrum (orthogonal basis for their
$L^2$ space consisting of group characters).
This disproves the Universal Spectrum Conjecture of Lagarias and Wang \cite{lagwang}.
Further, we construct a set in some finite Abelian group, which tiles the group but has no spectrum.
We extend this last example to the groups $\ZZ^d$ and $\RR^d$ (for $d \ge 5$)
thus disproving one direction of the Spectral Set Conjecture of Fuglede \cite{fuglede}.
The other direction was recently disproved by Tao \cite{tao}. 
\end{abstract}

\maketitle



\section{Introduction}\label{sec:intro}
Let $G$ be a locally compact Abelian group and $\Omega\subseteq G$ be a bounded open set.
We call $\Omega$ {\em spectral} if there is a set $\Lambda$ of continuous characters of
$G$ which forms an orthogonal basis for $L^2(\Omega)$. Such a set $\Lambda$ is called
a {\em spectrum} of $G$.
This paper is about a conjecture of Fuglede \cite{fuglede} (the {\em Spectral Set Conjecture}),
which states that a domain $\Omega$ in $\RR^d$ is spectral
if and only if it can tile $\RR^d$ by translation.
A set $\Omega$ tiles $\RR^d$ by translation
if there exists a set $T\subseteq\RR^d$ (called a {\em tiling complement} of $\Omega$) of translates
such that $\sum_{t\in T} \chi_\Omega(x-t) = 1$, for almost all $x\in\RR^d$.
Here $\chi_\Omega$ denotes the indicator function of $\Omega$.

Tao \cite{tao} has recently proved that the direction ``spectral $\Rightarrow$ tiling'' does not
hold (in dimension 5 and higher -- Matolcsi \cite{mat} has reduced this dimension to 4).
Here we prove that the direction ``tiling $\Rightarrow$ spectral'' is also false in dimension
5 and higher.

The Spectral Set Conjecture has attracted a lot of attention in the last decade.
Until Tao's example \cite{tao} there have been many results
\cite{fuglede,convexcurv,convexplane,spectralsym,konyagin,lagwang,lagszab,unispectra}
for special cases of domains, tiling complements
or spectra, all of them in the direction of supporting the conjecture.
For example, Fuglede showed that the conjecture is true if either the tiling
complement or the spectrum is assumed to be a lattice.
It may still be true that the conjecture is true for some rather large natural class
of domains, such as the convex domains \cite{convexplane}.
The reader should consult the references given in \cite{tao} or the survey \cite{milano}.

Our strategy is as follows.
The Spectral Set Conjecture makes sense in finite groups as well and we first disprove
the direction ``tiling $\Rightarrow$ spectral'' in
an appropriate finite group, just as Tao \cite{tao} did with the other direction.
This we do in \S \ref{sec3}, by first finding a counterexample to the
Universal Spectrum Conjecture of Lagarias and Wang \cite{lagwang}.
(This conjecture states that, in a finite group, if a set $T$ can tile the group
with tiling complements $T_1, \ldots, T_n$ these sets are all spectral and share
a common spectrum.)
We then produce in \S \ref{sec4}, using the example in the finite group,
a counterexample in the group $\ZZ^d$ and finally in $\RR^d$, where the Spectral
Set Conjecture was originally stated.

In \S \ref{sec2} we give necessary background material and describe notation.

\section{Preliminaries}\label{sec2}
Suppose $\Omega$ is a bounded open set in a locally compact Abelian group $G$.
We will only be interested in finite groups, $\ZZ^d$ and $\RR^d$ and everything we say from
now on applies to them.

We call $\Omega$ {\em spectral} if $L^2(\Omega)$
has an orthonormal basis
$$
\Lambda \subseteq \ft{G}
$$
of characters ($\ft{G}$ denotes the dual group of $G$ \cite{rudin}).
The set $\Lambda$ is then called a {\em spectrum} for $\Omega$.
In the groups we care about the characters are functions of the type
$x \to \exp(2\pi i \inner{\nu}{x})$, where $\nu$ takes values in an appropriate
subgroup of the torus $\TT^d$ (if $G$ is discrete) or in Euclidean space.

The inner product and norm on $L^2(\Omega)$ are
$$
\inner{f}{g}_\Omega = \int_\Omega f \overline{g},
\ \mbox{ and }\
\Norm{f}_\Omega^2 = \int_\Omega \Abs{f}^2.
$$
If $\lambda, \nu \in \ft{G}$ we have
$$
\inner{\lambda}{\nu}_\Omega = \ft{\chi_\Omega}(\nu-\lambda).
$$
which gives
$$
\mbox{$\Lambda$ is an orthogonal set} \Leftrightarrow
 \forall \lambda,\mu\in\Lambda, \lambda\neq\mu:
   \ \ft{\chi_\Omega}(\lambda-\mu) = 0
$$
For $\Lambda$ to be complete as well we must in addition have (Parseval)
\beql{parseval}
\forall f\in L^2(\Omega):\ \ 
    \Norm{f}_2^2 =  \frac{1}{\Abs{\Omega}^2} \sum_{\lambda\in\Lambda} \Abs{\inner{f}{\lambda}}^2.
\eeq
For the groups we care about (finite groups, $\ZZ^d$ and $\RR^d$)
in order for $\Lambda$ to be complete
it is sufficient to have \eqref{parseval} for any character $f \in \ft{G}$,
since then we have it in the closed linear span of these functions, which is all
of $L^2(\Omega)$.
An equivalent reformulation for $\Lambda$ to be a spectrum of $\Omega$
is therefore that
\beql{power-spectrum-tiling}
\sum_{\lambda\in\Lambda} \Abs{\ft{\chi_\Omega}}^2(x-\lambda) = \Abs{\Omega}^2,
\eeq
for almost every $x\in \ft{G}$.
For finite sets $\Omega$ (the group is finite or $\ZZ^d$) for a set $\Lambda \subseteq \ft{G}$
to be a spectrum it must satisfy the two conditions:
\begin{itemize}
\item[(a)] $\Lambda-\Lambda \subseteq \Set{\ft{\chi_\Omega}=0} \cup \Set{0}$, and
\item[(b)] $\# \Lambda = \# \Omega$.
\end{itemize}
For subsets $\Omega \subseteq \RR^d$, when the spectra are infinite, we fall back
on \eqref{power-spectrum-tiling}.

If $f \ge 0$ is in $L^1(G)$ and $T \subseteq G$ we say that $f$ tiles with $T$ at level
$\ell$ if $\sum_{t \in T} f(x-t) = \ell$ for almost all $x\in G$.
We denote this by ``$f+T = \ell G$'' and we call $T$ a {\em tiling complement} of $f$.
If $f = \chi_\Omega$ is the indicator function of some set then we just write
$\Omega+T = \ell G$ instead of $\chi_\Omega + T = \ell G$, and, in this case, if
$\ell$ is not specified it assumed to be 1.

In the finite group case it is immediate to show that $f+T$ is a tiling of $G$
if and only if
\beql{tiling-in-fourier}
\Set{\ft{f} = 0} \cup \Set{\ft{\chi_T} = 0} \cup \Set{0} = \ft{G}.
\eeq
There are analogs of this relationship that hold in the infinite case as well but we will not need
these here (see \cite{milano}).

If $f$ is a continuous function we write $Z(f)$ for its zero set. For a set $A$
we write $Z_A$ for the zero set of the Fourier Transform of its indicator function
$Z(\ft{\chi_A})$.

We now describe a generalization of the composition construction
appearing in \cite{mat}, Proposition 2.1. 

\begin{proposition}\label{prop1}
Let $G$ be a finite Abelian group, and $H\le G$ a subgroup.
Let $T_1, T_2, \dots T_k\subset H$ be subsets of $H$ such that they share a common tiling set in $H$;
i.e. there exists a set $T'\subset H$ such that $T_j+T'=H$ is a tiling for all $1\le j\le k$.
Consider any tiling decomposition $S+S'=G/H$ of the factor group $G/H$, with $\# S=k$, and take arbitrary representatives $s_1, s_2, \dots s_k$ from the cosets of $H$ corresponding to the set $S$. Then the set $\C :=\cup_{j=0}^{k}(s_j+T_j)$ is a tile in the group $G$.
\end{proposition} 
\begin{proof}
The proof is simply the observation that for any system of representatives $\tilde{S'}:=\{ s'_1, s'_2, \dots \}$ of $S'$ the set $T'+\tilde{S'}$ is a tiling set for $\C$ in 
$G$. 
\end{proof}
Despite the proof being obvious, this construction seems to include a large class of tilings
and it leads to some interesting examples.

A similar construction applies to spectral sets, as well.

\begin{proposition}\label{prop2}
Let $G$ be a finite Abelian group, and $H\le G$ a subgroup. Let $T_1, T_2, \dots T_k\subset H$ be subsets of $H$ such that they share a common spectrum in $H$; i.e. there exists a set $L\subset H$ such that $L$ is a spectrum of $T_j$ for all $1\le j\le k$.
Consider any spectral pair $(Q, Q')$ in the factor group $G/H$, with $\# Q=k$, and take arbitrary representatives $q_1, q_2, \dots q_k$ from the cosets of $H$ corresponding to the set $q$. Then the set $\C :=\cup_{j=0}^{k}(q_j+T_j)$ is spectral in the group $G$. 
\end{proposition}
\begin{proof}
We do not give a detailed proof of this statement, as we will not directly use it in the forthcoming arguments. Let us mention only that a spectrum of $\C$ is the set $L+\tilde{Q'}$, where $\tilde{Q'}$ denotes a system of representatives of $Q'$. The detailed proof proceeds along the same lines as in \cite{mat}, Proposition 2.1. 
\end{proof}

The main point of the two preceding constructions is that they are not entirely ``compatible''.
That is, one can hope to find sets $T_1, \dots T_k\subset H$
sharing a common tiling complement $T'$ but not sharing a common
spectrum $L$.
This would be a counterexample to the Universal Spectrum Conjecture.
Then the construction of Proposition \ref{prop1} will lead to a set $\C$ which tiles $G$,
but there is nothing to guarantee that $\C$ is spectral in $G$
(in fact, we will find a way to guarantee that $\C$ is {\em not} spectral).
This is exactly the route we will follow in \S \ref{sec3}. 

\section{Counterexamples in finite groups}\label{sec3}
Here we follow the path outlined in \S \ref{sec2} in order to produce an
example of a set $\C$ in a finite group $G$, such that $\C$ is a tile but is not spectral in $G$. 

The first step is to find a counterexample to the Universal Spectrum Conjecture.
We are looking for a finite group $G$ and a tile $T'$ in $G$ such that the tiling
complements $T_1, \dots T_k$ of $T'$ do not posess a common spectrum $L$. 

For a given $T'\subset G$, one sufficient condition for
the existence of a universal spectrum $L$, as pointed out in \cite{lagszab},
is to ensure that 
\begin{equation}\label{lagarias}
\# L \cdot \# T' = \# G\ \ \mathrm{and} \ \ L-L\subset Z^c_{T'}.
\end{equation}
Indeed, any tiling complement $T_j$ of $T'$ must satisfy $Z_{T_j}\supset  Z^c_{T'}\setminus \{ 0\}$, therefore condition \eqref{lagarias} ensures that $L$ is a spectrum of $T_j$.
(We do not know whether condition \eqref{lagarias} is also necessary
for the existence of a universal spectrum,
as suggested in \cite{lagszab} in the remarks following Theorem 3.1.)
Therefore, when trying to construct a set $T'$ having no universal spectrum,
one must exclude the existence of a set $L$ satisfying \eqref{lagarias}. 

Notice, that if $L$ satisfies \eqref{lagarias} then $L$ is not only a universal spectrum
for all tiling complements of $T'$, but also a universal tiling set of all spectra of $T'$.
Indeed, for any spectrum $Q$ of $T'$ we have $Q-Q\subset  Z_{T'}\cup \{ 0 \}$
therefore $L-L\cap Q-Q=\{ 0\}$ and $\# L\cdot \# Q$=$\#G$,
which ensures that $L+Q=G$ is a tiling. 
  
Having this observation in mind, one way to exclude the possibility of \eqref{lagarias}
is to choose a set $T'$ which posseses a particular spectrum $Q$
which does not tile $G$
(but recall that $T'$ itself must tile $G$ otherwise
the notion of universal spectrum is meaningless).
In other words, in some group $G$ take a spectral set $Q$ which does not tile $G$
(such examples already exist, cf. \cite{tao}, \cite{mat})
and choose any spectrum of $Q$ as a candidate for $T'$.
However, the examples in \cite{tao} and \cite{mat} are such that
$\# Q$ does not divide $\# G$,
therefore any choice for $T'$ is also doomed by divisibility reasons,
because $T'$ cannot tile $G$ either.
We circumvent this problem by increasing the size of the group $G$. 

The ideas above are summarized in the following proposition, which disproves
the Universal Spectrum Conjecture.
\begin{proposition}\label{usc}
Consider  $G=\Z_{6}^5$ and $E =\{ \vec{0}, \vec{e_1}, \vec{e_2}, \dots \vec{e_5}\}$ where $\vec{e_j}=(0, \dots 1, \dots , 0)^\top$.
The set $E$ tiles $G$ but has no universal spectrum in $G$. 
\end{proposition}
\begin{proof}
 
The existence of a universal spectrum $L$ is equivalent to the conditions $\# L=6^4$ and $L-L\subset (\bigcap_{j}Z_{T_j})\cup \{0\}$, where $T_j$ are all the tiling complements of $E$. 

Consider the set $K\subset G$ (from \cite[Theorem 3.1]{mat})
consisting of the columns of the following matrix
\beql{matrix-K}
K:=\left (
\begin{array}{cccccc}
0&0&2&2&4&4\\
0&2&0&4&4&2\\
0&2&4&0&2&4\\
0&4&4&2&0&2\\
0&4&2&4&2&0
\end{array}
\right )
\eeq
(We remark, that $K$ is a spectrum of $E$. In fact, this follows from the fact that the matrix 
\beql{matrix-K-prime}
K':= \frac{1}{3}\left (
\begin{array}{cccccc}
0&0&0&0&0&0\\
0&0&1&1&2&2\\
0&1&0&2&2&1\\
0&1&2&0&1&2\\
0&2&2&1&0&1\\
0&2&1&2&1&0
\end{array}
\right )
\eeq
is log-Hadamard.
This is a matrix $M$ for which the matrix $U_{ij} = \exp(2\pi i M_{ij})$ is orthogonal.
We will not use the fact that $K$ is a spectrum,
but it reflects the considerations preceding the proposition.)

Observe that $K$ is contained in the subgroup $H$ of vectors having even coordinates only.
However, $\# H=3^5$ and $\# K=6$, therefore $K$ cannot tile $H$ and, consequently,
it cannot tile $G$ either.
(It is easy to see that if a set tiles a group then it tiles the subgroup it
generates.)
It is also easy to check that the set $K-K$ consists of $\vec{0}$ and
only coordinate permutations of the vector $(0,2,2,4,4)^\top$.    
(In fact $K-K$ contains {\em all} coordinate permutations of $(0,2,2,4,4)^\top$,
but we do not need this.)

Next we show that $E$ admits some tiling complements $T_0, \dots T_{14}$,
which have no common spectrum.

Take the vector $\vec{v_1}=(1,2,3,4,5)^\top$
and define a group homomorphism $\f:G\to \Z_6$ by
$$
\f (\vec{x}):=\vec{v_1^\top}\vec{x}\ (\bmod\ 6).
$$
Then $\f$ is one-to-one on $E$, and the image of $E$ is the whole group $\Z_6$.
Therefore $T_0 = \ker\f$ is a tiling complement for $E$.
Notice that $Z_{T_0}^c$ contains all multiples of $\vec{v_1}$, and, in particular,
it contains
$2\vec{v_1}=(2,4,0,2,4)^\top$, and $4\vec{v_1}=-2 \vec{v_1} = (4,2,0,4,2)^\top$.
By appropriate permutations of the coordinates of $\vec{v_1}$
we can define vectors $\vec{v_2}, \dots ,\vec{v_{15}}$
and corresponding tiling sets $T_0, \dots ,\T_{14}$
in such a way that $( \bigcap_{j=0}^{14}Z_{T_j})\cap (K-K)=\{ 0\}$.
Therefore, a set $L$ satisfying  $\# L=6^4$ and
$L-L\subset (\bigcap_{j=0}^{14}Z_{T_j})\cup \{0\}$ cannot exist because
in that case $L+K$ would be a tiling of $G$, and we already know that $K$ is not a tile.
\end{proof}

Having found a counterexample to the Universal Spectrum Conjecture,
we use the construction of \S \ref{sec2} to exhibit the
failure of the Spectral Set Conjecture in finite groups.

\begin{proposition}\label{fug}
Consider  $G_2=\Z_{6}^5\times \Z_{15}$ and
$\C=\bigcup_{j=0}^{14}(f_j+\tilde{T_j})$, where $f_j=(0,0,0,0,0,j)^\top$
and $\tilde{T_j}$ are the sets appearing in
Proposition \ref{usc} extended by 0 as the last coordinate.
Then $\C$ is a tile in $G_2$ but it is not spectral. 
\end{proposition} 
\begin{proof}
In this proof the notation $\tilde{A}$ always refers to a set $A\subset G=\Z_6^5$ extended by 0 as the last coordinate.

The fact that $\C$ is a tile follows from Proposition \ref{prop1}
or can easily be seen directly: the tiling complement of $\C$ is $\tilde{E}$. 

To see that $\C$ is not spectral, note first that the set $\tilde{K}$
is contained in the subgroup $\tilde{H}$ (defined in the proof of Proposition \ref{usc}),
therefore it cannot tile $G_2$ because of divisibility reasons.

Any spectrum $Q$ of $\C$ must satisfy $\# Q = \# \C =6^4\cdot 15$
and $Q-Q\subset Z_\C\cup \{ 0\}$.
Consider the vector ${\vec{\tilde{k_1}}}=(0,2,2,4,4,0)^\top\in \tilde{K}-\tilde{K}$.
We show that $\vec{\tilde{k_1}}\notin Z_{\C}$.
Indeed,
$$
\ft{\chi}_{\C}({\vec{\tilde{k_1}}})= \sum_{j=0}^{14}\ft{\chi}_{T_j}(\vec{k_1})>0
$$
because each term is nonnegative (each $T_j$ being a subgroup in $G$),
and at least one term is non-zero by construction.
The same argument shows that  $\vec{\tilde{k_j}}\notin
Z_{\C}$ for all $\vec{k_j}\in K-K$. 

Therefore, any spectrum $Q$ of $\C$ must satisfy
$(Q-Q)\cap (\tilde{K}-\tilde{K})=\{ 0\}$,
which, since $\#Q \cdot \#\tilde{K} = \#G_2$, implies
$Q+\tilde{K}=G_2$.
This is a contradiction since $\tilde{K}$ is not tile in $G_2$   
\end{proof}

\section{Transition to $\Z^d$ and $\R^d$}\label{sec4}
We now describe a general transition scheme from the finite group setting
to $\Z^d$ and $\R^d$.
As a result we find a set in $\R^6$, which is a finite union of unit cubes
(placed at points with integer coordinates),
which tiles $\R^6$ by translations but is not spectral.    

First we prove this in the group $\ZZ^d$.
\begin{theorem}\label{th:zd}
Let $\vec n = (n_1,\ldots,n_d) \in \ZZ^d$ and suppose that the set
$A \subset G = \ZZ_{n_1}\times\cdots\times\ZZ_{n_d}$ is not spectral
in the group $G$.
Write
\beql{ft-grid}
T = T(\vec n, k) = \Set{0,n_1,2n_1,\ldots,(k-1)n_1}\times\cdots\times
			\Set{0,n_d,2n_d,\ldots,(k-1)n_d},
\eeq
and define $A(k) = A + T$.
Then, if $k$ is large enough,  the set $A(k)\subset\ZZ^d$ is not spectral in $\ZZ^d$.
\end{theorem}

\begin{proof}
We will argue by contradiction and assume that there are arbitrarily large
values of $k$ for which $A(k)$ is spectral in $\ZZ^d$.

Observe first that $\chi_{A(k)} = \chi_A * \chi_T$, hence,
writing $Z(f) = \Set{f = 0}$, we obtain
$$
Z(\ft{\chi_{A(k)}}) = Z(\ft{\chi_A}) \cup Z(\ft{\chi_T}).
$$
Elementary calculation of $\ft{\chi_T}$ (it is a cartesian product) shows that
it is a union of ``hyperplanes''
\beql{zt}
Z(\ft{\chi_T}) =
 \Set{\vec \xi \in \TT^d:\ \exists j\ \exists \nu\in\ZZ,\mbox{ $k$ does not divide $\nu$, such that } \xi_j = \frac{\nu}{kn_j}}.
\eeq
Define the group
$$
H = \Set{\vec \xi \in \TT^d:\ \forall j\ \exists \nu\in\ZZ\ \mbox{such that }
	\xi_j = \frac{\nu}{n_j}}.
$$
which is the group of characters of the group $G$ and does not depend on $k$.
Observe that $H+(Q-Q)$ does not intersect $Z(\ft{\chi_T})$, where
$$
Q = \left[0,\frac{1}{kn_1}\right)\times\cdots\times\left[0,\frac{1}{kn_d}\right).
$$

Assume now that $S\subseteq\TT^d$ is a spectrum of $A(k)$, so that $\# S = \# A(k) = r k^d$,
if $r = \# A$.
Define, for $\vec \nu \in \Set{0,\ldots,k-1}^d$, the sets
$$
S_{\vec \nu} = \Set{\vec\xi \in S:\ \vec\xi \in (\frac{\nu_1}{k n_1},\ldots,\frac{\nu_d}{k n_d})
		+ Q + (\frac{m_1}{n_1},\ldots,\frac{m_d}{n_d}),\ \mbox{for some $\vec m \in \ZZ^d$}}.
$$
\comment{
It follows from \eqref{few-bad-zeros} that
$$
\#\left(S \setminus \bigcup_{\vec \nu} S_{\vec \nu} \right) \le C k^{d-1}.
$$
}
Since the number of the $S_{\vec \nu}$ is $k^d$ it follows that
there exists some $\vec \mu$ for which $\# S_{\vec \mu} \ge r$.

We also note that, if $k$ is sufficiently large, then any translate of $Q$ 
may contain at most one point of the spectrum. The reason is that $Q-Q$ contains
no point of $Z(\ft{\chi_T})$ (for any $k$) and no point of $Z(\ft{\chi_A})$
for all large $k$ (as $\ft{\chi_A}(\vec 0) > 0$).

Observe next that if $\vec x, \vec y \in S_{\vec \mu}$ then
\begin{eqnarray*}
\vec x - \vec y &\in& H + (Q-Q)\\
 &=& H + \left(-\frac{1}{kn_1},\frac{1}{kn_1}\right)\times\cdots\times\left(-\frac{1}{kn_d},\frac{1}{kn_d}\right)
\end{eqnarray*}
and that this set does not intersect $Z(\ft{\chi_T})$ (from \eqref{zt}). 
It follows that for all $\vec x, \vec y \in S_{\vec \mu}$ we have
$\vec x - \vec y \in Z(\ft{\chi_A})$. 

Let $k$ be sufficiently large so that for all points $\vec h \in H$ for
which $\ft{\chi_A}(\vec h) \neq 0$ the rectangle $\vec h + Q-Q$ does not
intersect $Z(\ft{\chi_A})$.
For each $\vec x \in \TT^d$ define $\lambda(\vec x)$ to
be the unique point $\vec z$ whose $j$-th coordinate is an integer multiple
of $\frac{1}{k n_j}$ for which $\vec x \in \vec z + Q$.
If $\vec x, \vec y \in S_{\vec \mu}$ it follows that
$\lambda(\vec x)-\lambda(\vec y) \in H \cap Z(\ft{\chi_A})$
(otherwise $\lambda(\vec x)-\lambda(\vec y) + Q-Q$ would contain
no zeros of $\ft{\chi_A}$ and $\vec x - \vec y$ would not
be in $Z(\ft{\chi_A})$).
Define now $\Lambda = \Set{\lambda(\vec x):\ \vec x \in S_{\vec \mu}}$.
It is obvious that $\# \Lambda \ge r$ and
$\Lambda-\Lambda \subseteq Z(\ft{\chi_A}) \cup \Set{\vec 0}$,
therefore $\Lambda$ is a spectrum of $A$, which is a contradiction.

\comment{
Now let $k \to \infty$ and pass to a subsequence to ensure that the sets
$S_{\vec \mu}$, properly translated (for example, translate them to contain $\vec 0 \in \TT^d$),
converge to a subset $\Lambda$ of $H$.
It follows that $\# \Lambda \ge r$ and $\Lambda-\Lambda \subseteq Z(\ft{\chi_A})$,
by the continuity of $\ft{\chi_A}$.
This implies that $\Lambda$ is a spectrum of $A$ (and that $\# \Lambda = r$),
a contradiction.
}
\end{proof}

Non-spectral tiles can be pulled from $\ZZ^d$ to $\RR^d$ using the following.
\begin{theorem}\label{th:discrete-continuous}
Suppose $A \subseteq \ZZ^d$ is a finite set and $Q = [0,1)^d$. Then
$A$ is a spectral set in $\ZZ^d$ if and only if $A+Q$ is a spectral set in $\RR^d$.
\end{theorem}

\begin{proof}
Write $E=A+Q$. Then $\ft{\chi_E} = \ft{\chi_A} \ft{\chi_Q}$ and
$Z(\ft{\chi_E}) = Z(\ft{\chi_A}) \cup Z(\ft{\chi_Q})$.
By calculation we have
$$
Z(\ft{\chi_Q}) =
	\Set{\vec \xi \in \RR^d:\ \exists j \ \mbox{such that } \xi_j \in \ZZ\setminus\Set{\vec 0}}.
$$

Now suppose $\Lambda \subset \TT^d$ is a spectrum of $A$ as a subset of $\ZZ^d$.
Viewing $\TT^d$ as $Q$ we observe that the set $Z(\ft{\chi_A})$ is periodic with
$\ZZ^d$ as a period lattice.
Define now $S = \Lambda + \ZZ^d$.
The differences of $S$ are either points which are on $Z(\ft{\chi_A})$ (mod $\ZZ^d$)
or points with all integer coordinates. In any case these differences fall in 
$Z(\ft{\chi_E})$, hence $\Abs{\ft{\chi_E}}^2+S$ is a packing of $\RR^d$ at level $(\# A)^2$.
Furthermore, the density of $S$ is $\# A$ which, along with the periodicity of $S$, implies
that $\Abs{\ft{\chi_E}}^2+S$ is a tiling of $\RR^d$ at level $(\# A)^2$.
That is, $S$ is a spectrum for $E$.

Conversely, assume $S$ is a spectrum for $E$ as a subset of $\RR^d$.
It follows that the density of $S$ is equal to $\Abs{E} = \# A$, hence there exists $\vec k \in \ZZ^d$
such that $\vec k + Q$ contains at least $\# A$ points of $S$.
Call the set of these points $S_1$, and observe that the differences of points of $S_1$
are contained in $Q-Q = (-1,1)^d$, and that $Q-Q$ does not intersect $Z(\ft{\chi_Q})$.
It follows that the differences of the points of $S_1$ are all in $Z(\ft{\chi_A})$,
and, since their number is $\# A$, they form a spectrum of $A$ as a subset of $\ZZ^d$.
\end{proof}

In conclusion we have proved the following.
\begin{theorem}\label{th:last}
In each of the groups $\ZZ^d$ and $\RR^d$, $d\ge 5$, there exists a set which
tiles the group by translation but is not spectral
\end{theorem}

\begin{proof}
It is easy to see that if $A$ is our non-spectral tile in the finite group
$\ZZ_6^5\times\ZZ_{15}$ then the set $A(k) \subseteq \ZZ^6$ which appears in Theorem \ref{th:zd}
is a tile as well, and by that Theorem it is not spectral.
Using Theorem \ref{th:discrete-continuous} we can construct a set with
these properties in $\RR^6$ by adding a unit cube at each point.

To get down to dimension 5, notice that the construction in Theorem \ref{fug}
can be repeated verbatim in the group $G_3 = \ZZ_6^5 \times \ZZ_{17}$
instead of the group $\ZZ_6^5 \times \ZZ_{15}$. (Just repeat the set $T_{15}$ two more times.)
But now we view $G_3$ as the group $\ZZ_6^4 \times \ZZ_{6\cdot17}$ (as 6 and 17 are coprime).
Theorems \ref{th:zd} and \ref{th:discrete-continuous} now give examples in dimension 5.

If $d\ge 6$ then the set of $\ZZ^6$ which is tiling but not spectral will still
be such in $\ZZ^d$ when viewed as a subset of that in the obvious way. The preservation
of tiling property is obvious, and one can easily show that the existence of any spectrum
in $\TT^d$ implies the existence of a spectrum in $\TT^6$. 
We go to $\RR^d$ again using Theorem \ref{th:discrete-continuous}.
\end{proof}

\end{document}